\documentclass[runningheads, 12pt]{llncs}

\usepackage{amssymb}
\usepackage [english]{babel}
\usepackage[utf8]{inputenc}
\usepackage[T1]{fontenc}

\usepackage{amssymb}
\usepackage{amsmath}
\usepackage{amsfonts}

\usepackage{graphicx}
\usepackage{url}

\usepackage[usenames]{color}

\usepackage{algorithm}
\usepackage{algpseudocode}
\usepackage{parskip}
\usepackage{subcaption}
\usepackage{makecell}
\usepackage{xfrac}

  \usepackage{geometry}
\newcommand{\PathToPic}{}

 \newcommand{\MatTwo}[4]
  {\begin{bmatrix}
    #1 & #2\\
    #3 & #4
  \end{bmatrix}}
  
  \newcommand{\VecTwo}[2]
  {\begin{bmatrix}
    #1 \\
    #2 
  \end{bmatrix}}

    \newcommand{\keywords}[1]{\par\addvspace\baselineskip  
     \noindent\keywordname\enspace\ignorespaces#1}

\begin{document}

\urldef{\mailsa}\path|christoph.hofer@jku.at|
\urldef{\mailsb}\path|ulrich.langer@jku.at|
\urldef{\mailsc}\path|stefan.takacs@ricam.oeaw.ac.at|

\title{Inexact Dual-Primal Isogeometric Tearing and Interconnecting Methods}
\author{Christoph Hofer$^1$, Ulrich Langer$^{1,2}$ and Stefan Takacs$^2$}
\institute{ $^1$ Johannes Kepler University (JKU),
				Altenbergerstr. 69, A-4040 Linz, Austria,\\
christoph.hofer@jku.at, ulrich.langer@jku.at\\
$^2$ Austrian Academy of Sciences, RICAM,
Altenbergerstr. 69, A-4040 Linz, Austria,\\
ulrich.langer@ricam.oeaw.ac.at, stefan.takacs@ricam.oeaw.ac.at
 }

\noindent
\maketitle

\begin{abstract}
In this paper, we investigate inexact variants of dual-primal isogeometric tearing and interconnecting methods for solving large-scale systems of linear equations arising from Galerkin isogeometric discretizations of elliptic boundary value problems.
The considered methods are extensions of standard finite element tearing and interconnecting methods to isogeometric analysis. The algorithms are implemented by means of energy minimizing primal subspaces. 
We discuss the replacement of local  sparse direct solvers by iterative methods, particularly, multigrid solvers. 
We investigate the incorporation of 
 these
iterative solvers into different 
formulations 
of the algorithm. Finally, we present numerical examples comparing 
 the performance of these 
inexact versions.
\end{abstract}


\pagestyle{myheadings}
\thispagestyle{plain}
\markboth{}{C. Hofer, U. Langer and S. Takacs, Inexact IETI-DP}

\section{Introduction}
Isogeometric Analysis (IgA) is a novel methodology for the numerical solution of partial differential equations (PDEs). IgA was first introduced by Hughes, Cottrell and Bazilevs in \cite{HL:HughesCottrellBazilevs:2005a}, see also
the survey article \cite{HL:BeiraodaVeigaBuffaSangalliVazquez:2014a}. 
In IgA, 
for both  the representation of the geometry and the approximation of the solution,
spline-based spaces are chosen. The most common choices are B-Splines, Non Uniform Rational B-Splines (NURBS), T-Splines, Truncated Hierarchical B-Splines (THB-Splines), see \cite{HL:BeiraodaVeigaBuffaSangalliVazquez:2014a} and references therein.
One of the strengths of IgA 
 consists in
its capability of creating high-order 
 smooth
function spaces, while keeping the number of degrees of freedom small. Originally, IgA was formulated 
 by means of
one \emph{global geometry mapping}, which restricts the method to simple domains being topologically equivalent to the unit square or the unit cube. More complicated domains are represented by decomposing them into such simple domains, called \emph{patches} or \emph{subdomains}. In such a \emph{multi-patch} setting,
each of the patches has its own geometry mapping, and all of the patches can be discretized
 by the use of
different spline spaces.
 
 We are interested in fast solvers for linear systems arising from the discretization of elliptic PDEs by means of IgA.
 We investigate non-overlapping domain decomposition (DD) methods of the dual-primal tearing and interconnecting type. These methods are closely related to the Balancing Domain Decomposition by Constraints (BDDC) methods, see \cite{HL:ToselliWidlund:2005a,HL:Pechstein:2013a} and references therein. The version based on a conforming Galerkin (cG) discretization, called  dual-primal isogeometric tearing and interconnecting (IETI-DP) method, was first introduced in \cite{HL:KleissPechsteinJuettlerTomar:2012a}. 
 The related IgA BDDC method was analyzed in \cite{HL:VeigaChoPavarinoScacchi:2013a}. Typically, the local problems are solved using a sparse Cholesky factorization.
 However, especially in IgA, one may run out of memory for big problems. A remedy would be to use inexact solvers for the local subproblems, as introduced in \cite{HL:KlawonnRheinbach:2007a}. The aim of this work is to investigate how 
 local sparse direct solvers can be replaced by 
inexact methods, like multigrid (MG). This leads to several different variants of the IETI-DP algorithm, each with its own advantages and disadvantages.

In the present paper, we consider the following  weak formulation of a
second-order elliptic boundary value problem (BVP)
in a bounded Lipschitz domain $\Omega\subset \mathbb{R}^d,$ with $d\in\{2,3\}$,
 as model problem:
Find $u \in V_0:=H^1_0(\Omega)$
such that
\begin{align}
  \label{equ:ModelVar}
    a(u,v) = \left\langle F, v \right\rangle \quad \forall v \in V_{0}.
\end{align}
The bilinear form
$a(\cdot,\cdot): V_0 \times V_{0} \rightarrow \mathbb{R}$
and the linear form $\left\langle F, \cdot \right\rangle: V_{0} \rightarrow \mathbb{R}$
are given by
\begin{equation*}
a(u,v) := \int_\Omega \nabla u \cdot \nabla v \,dx
\quad \mbox{and} \quad
\left\langle F, v \right\rangle := \int_\Omega f v \,dx,
\end{equation*}
 respectively. We assume that the given right hand side function $f$ is sufficiently smooth.
 
\section{Isogeometric Analysis and IETI-DP}

On the unit interval, 
for any spline degree $p$ and number of basis functions $M$,
 we define
the one dimensional B-Spline basis $(\widehat{N}_{i,p})_{i=1}^M$ via the Cox-De Boor's algorithm, cf. \cite{HL:BeiraodaVeigaBuffaSangalliVazquez:2014a}.
On the \emph{parameter domain} $\widehat{\Omega}:=(0,1)^d$,
a  multivariate basis is realized by the tensor product of such  univariate bases functions,
again denoted by $\widehat{N}_{i,p}$, where $i=(i_1,\ldots,i_d)\in\mathcal{I}:= \{1,\ldots,M_1\}\times\ldots\times\{1,\ldots,M_d\}$
and $p=(p_1,\ldots,p_d)$ are multi-indices.

In standard (single-patch) IgA, the \emph{physical domain} $\Omega$ is given as the image of the parameter domain under the
\emph{geometrical mapping} $G :\; \widehat{\Omega} \rightarrow \mathbb{R}^{{d}}$, defined by
	$    G(\xi) := \sum_{i\in \mathcal{I}} P_i \widehat{N}_{i,p}(\xi),$
with the control points $P_i \in \mathbb{R}^{{d}}$, $i\in \mathcal{I}$. 
In a multi-patch setting, the domain $\Omega$ (\emph{multipatch domain}) is decomposed into non-overlapping patches $\Omega^{(k)}$,  $k=1,\ldots,N$, 
such that $\overline{\Omega}:=\bigcup_{k=1}^N\overline{\Omega}^{(k)}$. Each patch $\Omega^{(k)}:=G^{(k)}(\widehat{\Omega})$ is represented 
 by
its own geometrical mapping. 
We call $\Gamma:=\bigcup_{k>l}\partial \Omega^{(k)}\cap \partial\Omega^{(l)}$ the \emph{interface},
and denote its restriction to one of the patches $\Omega^{(k)}$ by $\Gamma^{(k)}:= \Gamma\cap \partial \Omega^{(k)}$.
Here and in what follows, the superscript $(k)$ denotes the restriction of the underlying symbol to the patch $\Omega^{(k)}$.

We use the B-Splines not only for defining the geometry, but also for representing the approximate solution of the 
 BVP.
 Once
the basis functions are defined
on the parameter domain $\widehat{\Omega}$, we define the bases on the physical domain $\Omega^{(k)}$ via the standard pull-back principle, and obtain the basis functions
${N}_{i,p}:=\widehat{N}_{i,p}~\circ~G^{-1}$.
 
 The main idea of IETI-DP is to decouple the patches 
 by tearing the interface unknowns
  which introduces additional degrees of freedom (dofs). We denote the resulting space by
 $V_h$. Then, continuity is again enforced using Lagrange multipliers $\lambda$. Doing so, the local subproblems
 on each patch are essentially pure Neumann problems (at least for interior patches).
 Therefore, they
 have a kernel 
  consisting of the constant functions in our case.
 So, a Schur complement formulation is not possible. In order
 to overcome this problem, certain continuity conditions are enforced strongly, i.e., by incorporating into the space $V_h$, (\emph{strongly enforced continuity conditions}) which yields the smaller space $\widetilde{V}_h$.
 There, we formulate the following problem.
 Find $(u,\lambda) \in \widetilde{V}_h\times \Lambda$ such that
       \begin {align}
       \label{equ:saddlePoint}
      	\MatTwo{\widetilde{K}}{\widetilde{B}^T}{\widetilde{B}}{0} \VecTwo{u}{\lambda} = \VecTwo{\widetilde{f}}{0},
      \end {align}
 where $\widetilde{K}$ is the stiffness matrix, $\widetilde{B}$ the jump operator, and $\widetilde{f}$ the right hand side, all in $\widetilde{V}_h$.
 
 As next step, we split $V_h$ into interior dofs and interface dofs, which yields an interface space $W$. By splitting $\widetilde{V}_h$ analogously, we
 obtain the space $\widetilde{W}$. Based on this splitting, we
 formulate the problem using the Schur complement of the stiffness matrix $K$ in $V_h$ with respect to the interface dofs:
 $S:= K_{BB}-K_{BI}(K_{II})^{-1}K_{IB}$, where the subindices $B$ and $I$ denote the boundary and
 interior dofs, respectively. The restriction of $S$ into $\widetilde{W}$ is denoted by $\widetilde{S}$, which yields the saddle-point formulation
 of the problem: Find $(w,\lambda) \in \widetilde{W}\times \Lambda$ such that
       \begin {align}\label{equ:saddlePoint2}
        \MatTwo{\widetilde{S}}{\widetilde{B}^T}{\widetilde{B}}{0} \VecTwo{w}{\lambda} = \VecTwo{\widetilde{g}}{0},
      \end {align}
 where $\widetilde{g} := \widetilde{I}^T(f_B -  K_{BI} (K_{II})^{-1}f_I)$ and $\widetilde{I}$ 
 is the canonical embedding of $\widetilde{W}$ in $W$.
 
 We denote the subspace of $\widetilde{W}$ satisfying the strongly enforced continuity conditions homogeneously by $W_\Delta$ and
 the $S$-orthogonal complement by $W_\Pi$. In the literature, our choice of $W_\Pi$ is often called \emph{energy minimizing primal subspace}.
 Finally, we can define the Schur complement $F$ of the saddle-point problem~\eqref{equ:saddlePoint2}, and obtain the problem:
 Find $\lambda \in U$ such that
\begin{align}
   \label{equ:SchurFinal}
      F\lambda = d,
\end{align}
 where $F:= \widetilde{B} \widetilde{S}^{-1}\widetilde{B}^T$ and $d:= \widetilde{B}\widetilde{S}^{-1} \widetilde{g}$.

 Equation \eqref{equ:SchurFinal} is solved 
  by means of
 the conjugate gradient (CG) algorithm using the scaled Dirichlet preconditioner $M_{sD}^{-1} := B_DSB_D^T$, where $B_D$ is a scaled version of the jump operator $B$ on $V_h$. Note that we can approximate $\widetilde{S}^{-1}$ because
 $\widetilde{S}$ can be represented (by reordering of the dofs) as a block diagonal matrix, consisting of matrices $S_{\Delta\Delta}^{(k)}$ for each patch and the matrix
 $S_{\Pi\Pi}$. For a summary of the algorithm and a more detailed explanation, we refer, e.g., to \cite{HL:ToselliWidlund:2005a,HL:Pechstein:2013a,HL:HoferLanger:2016a} and references therein.

\section{Incorporating Multigrid in IETI-DP}
\label{sec:MGinIETI}

We investigate different possibilities 
to incorporate a multigrid solver into the IETI-DP algorithm. The application of the IETI-DP algorithm requires the solution of linear systems at certain places. Two types of local problems are involved: \emph{Dirichlet problems} and \emph{Neumann problems}.

\subsection{Local Dirichlet problems}
\label{sec:locD}
We have to solve linear systems with system matrix $K_{II}^{(k)}$ in the application of $S$ in the preconditioner and when
calculating the right hand side $\widetilde{g}$. These linear systems are Dirichlet problems. (They would have Neumann boundary conditions only
if the patch boundary contribute to the Neumann boundary of the whole domain.)
The right hand side $\widetilde{g}$ has to be computed very accurately, i.e., at least up to discretization error. 
However, for the preconditioner,
a few MG V-cycles are usually enough,
since we only have to ensure the spectral equivalence of the inexact scaled Dirichlet preconditioner to the exact one, cf.~\cite{HL:KlawonnLanserRheinbach:2016b}
and references therein,

\subsection{Local Neumann problems}
The second class of local problems are Neumann problems.
They appear in the construction of the $S$-orthogonal basis for $W_\Pi$ and in the application of $S_{\Delta\Delta}$.
Let us first investigate the construction of the basis $\{\phi^{(k)}_j\}_j$ for $W_\Pi^{(k)}$. Since we look for a nodal basis, which is $S$-orthogonal, we have to solve the following linear system
    \begin{align}\label{equ:saddlePointBasis:0}
     \MatTwo{S^{(k)}}{{C^{(k)}}^T}{C^{(k)}}{0}\VecTwo{\phi_j^{(k)}}{\mu_j^{(k)}} = \VecTwo{0}{\boldsymbol{e}_j^{(k)}}, \quad \forall j\in\{1,\ldots,n_{\Pi}^{(k)}\},
    \end{align}
    where $\boldsymbol{e}_j^{(k)} \in \mathbb{R}^{n_{\Pi}^{(k)}}$ is the $j$-th unit vector and the matrix $C^{(k)}$ realizes the $n_{\Pi}^{(k)}$ strongly enforced continuity conditions contributing to the patch $\Omega^{(k)}$. This system has to be solved for $n_{\Pi}^{(k)}$ right hand sides, which is an advantage for direct solvers over iterative solvers 
     because the expensive factorization must be computed only once.
    Instead of solving~\eqref{equ:saddlePointBasis:0}
    directly, we use the approach proposed in~\cite{HL:Pechstein:2013a}, solve
        \begin{align}
        \label{equ:saddlePointBasis}
     \MatTwo{K^{(k)}}{{{C}^{(k)}}^T}{{C}^{(k)}}{0}\VecTwo{\overline{\phi}_j^{(k)}}{\mu_j^{(k)}} = \VecTwo{0}{\boldsymbol{e}_j^{(k)}}, \quad \forall j\in\{1,\ldots,n_{\Pi}^{(k)}\},
    \end{align}
     and obtain the desired basis functions by $\phi_j=\overline{\phi}_{j}|_{\Gamma^{(k)}}$. Note that $\{\overline{\phi}^{(k)}_j\}_j$ is a $K$-orthogonal basis. 
    If the patch $\Omega^{(k)}$ does not touch the boundary $\partial\Omega$, the upper left block becomes semi-definite due to the presence of a kernel. We are looking for a way to use the CG algorithm. 
    As long as there is no kernel, i.e., where $\partial\Omega^{(k)}\cap\partial\Omega \neq \emptyset$, one straightforward way would be to use the Bramble-Pasciak conjugate gradient (BPCG) algorithm or one of its variations, see \cite{HL:BramblePasciak:1988a,HL:StollWathen:2008a}. However, these iterative methods require that the upper left block is positive definite. The remedy is a special preconditioner and a non standard inner product for the CG algorithm, leading to the \emph{Sch\"oberl-Zulehner} (SZ) preconditioner, see \cite{HL:SchoeberlZulehner:2007a}. An alternative approach  would be to use the MinRes method with a
    block diagonal preconditioner, for which our experiments indicated a larger number of iterations.
  
  The SZ preconditioner for~\eqref{equ:saddlePointBasis} requires preconditioners $\hat{K}^{(k)}$ and $\hat{H}^{(k)}$ for the upper left block $K^{(k)}$ and its inexact Schur complement $H^{(k)}:=C^{(k)}{(\hat{K}^{(k)})}^{-1} {C^{(k)}}^T$, respectively.
  The preconditioner $K^{(k)}$ shall be realized 
  by
  a few MG V-cycles. It is required that $\hat{K}^{(k)}>K^{(k)}$, which implies that $\hat{K}^{(k)}$ has to be positive definite. In order to handle also the case where $K^{(k)}$ is singular, we need to set up MG based on a regularized matrix
  $ K_M^{(k)}:= K^{(k)} + \alpha \widehat{M}^{(k)}, $
  where $\alpha$ is chosen to be $10^{-2}$ and $\widehat{M}^{(k)}$ is the mass matrix on the parameter domain. Note, we can exploit the tensor product structure to efficiently assemble the mass matrix $\widehat{M}^{(k)}$. Finally, this provides us with an appropriate preconditioner $\hat{K}^{(k)}$ for $K^{(k)}$.
  Secondly, the SZ preconditioner requires that $\hat{H}^{(k)}<H^{(k)}$. Since in our case the number of rows of $C^{(k)}$ is given by $n_{\Pi}^{(k)}$, a small
  number that does not change during refinement, we calculate the inexact Schur complement exactly. This can be performed by applying ${(\hat{K}^{(k)})}^{-1}$ 
   to  $n_{\Pi}^{(k)}$ vectors.
 Finally, by a suitable scaling, e.g., $\hat{H}^{(k)}:=0.99H^{(k)}$, we obtain the desired matrix inequality. Having the preconditioners $\hat{K}^{(k)}$ and $\hat{H}^{(k)}$, we apply CG with the
  SZ
  preconditioner to construct the basis for $W_\Pi^{(k)}$.
  
  The second type of Neumann problem appears in the application of $F$. We look for a solution of the system $S_{\Delta\Delta}^{(k)} w_\Delta^{(k)} = f_\Delta^{(k)}$, which can be written as
  \begin{align}
  \label{equ:localNeumann}
  	\MatTwo{S^{(k)}}{{C^{(k)}}^T}{C^{(k)}}{0}\VecTwo{w_\Delta^{(k)}}{\mu^{(k)}} = \VecTwo{f^{(k)}}{0}.
  \end{align}
Certainly, one can use the same method as above. However, we can utilize the fact that we search for a minimizer of
$\tfrac12 (S^{(k)} w^{(k)},w^{(k)}) - (w^{(k)},f^{(k)})$ in the subspace given by $C^{(k)} w^{(k)} = 0$. This solution can be computed by
first solving the unconstrained problem and projecting the minimizer into the subspace using a energy-minimizing projection.
The projection is trivial because the decomposition of $\widetilde{W}$ into $W_\Pi$ and $W_\Delta$ is $S$-orthogonal.

Note that the CG algorithm, when applied to a positive semidefinite matrix, stays in the factor space with respect to the kernel and
computes one of the minimizers. The solution of the constrained minimization problem is, as outlined above, obtained by applying the projection.
As long as the number of CG iterations is not too large, numerical instabilities are not observed when applying CG to a positive semidefinite problem.
  
  The $S$-orthogonal basis has to be computed very accurate in order to maintain the orthogonality. Because the equation
  $S_{\Delta\Delta}^{(k)} w_\Delta^{(k)} = f_\Delta^{(k)}$ appears in the system matrix $F$, its solution also requires an accuracy of at least the discretization error.

\subsection{Variants of inexact formulations}
\label{sec:formulations}

From the discussion above, we deduce four (reasonable) combinations of the IETI-DP method with direct solvers and MG.

\begin{description}
	\item (D-D) This is the classical IETI-DP method, where we use direct solvers everywhere. 
	\item (D-MG) We use MG in the scaled preconditioner for the solution of the local Dirichlet problems and the transformation of the right hand side, see Section~\ref{sec:locD}. As already mentioned, the required accuracy for computing $\widetilde{g}$ has to be of the order of discretization error, whereas for the preconditioner, 
	 a few
	V-cycles 
	are enough. 
	\item (MG-MG) We use MG for all patch local problems, i.e., the local Dirichlet and Neumann problems. This implies that also the calculation of the basis for $W_{\Delta}$ is performed by means of MG, which turns out to be very costly. Moreover, for each application of $F$, we have to solve a local Neumann problem in $W_{\Delta}$ with the accuracy in the order of the discretization error.
	\item (MG-MG-S) To overcome the 
	 efficiency
	problem of applying MG at each iteration up to a small precision, we use the saddle point formulation instead of $F$. On the one hand, at each iteration step we only have to apply a given matrix instead of solving a linear system. On the other hand, we  now have to deal 
	with a saddle point problem.
	 Moreover, the iteration is not only applied to the interface dofs, but also to the dofs in the whole domain.
\end{description}

 We will always assume that the considered multipatch domain has only a moderate number of patches, such that the coarse problem  can still be handled 
 by
 a direct solver. For extensions to inexact version for the coarse problem, we refer to, e.g., \cite{HL:KlawonnRheinbach:2007a}. 

      For the first three methods, we use the CG method to solve $F\lambda=d$ as outer iteration. For the last setting (MG-MG-S), we have to deal with the saddle point problem \eqref{equ:saddlePoint}, which we solve using the BPCG method. The building blocks for this method are a preconditioner $\hat{\widetilde{K}}$ for $\widetilde{K}$ and $\hat{F}$ for the Schur complement $F$. The construction of $\hat{\widetilde{K}}$ follows the same steps as in the previous section, but we only apply 
       a few
      MG V-cycles. Concerning $\hat{F}$, a good choice is the scaled Dirichlet preconditioner $M_{sD}^{-1}$, cf. \cite{HL:KlawonnRheinbach:2007a}. 

\section{Numerical Experiments}
\begin{figure}
      \begin{subfigure}{0.45\textwidth}     
        \center{\includegraphics[width=0.5\textwidth]{\PathToPic 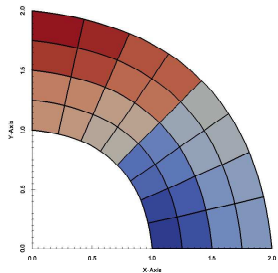}\caption{Quarter annulus}}
          
  \end{subfigure}
        \begin{subfigure}{0.45\textwidth}
        \center{\includegraphics[width=0.5\textwidth]{\PathToPic 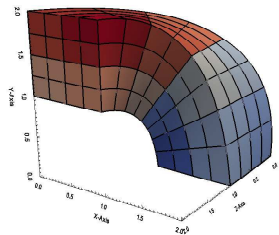}\caption{Twisted quarter annulus}}
  \end{subfigure}
  \caption{Illustration of the two and three dimensional computational domain.}
  \label{fig:QA}
\end{figure}

We solve the model problem \eqref{equ:ModelVar} on a two and a three dimensional computational domain.
In the two dimensional case, we use the quarter annulus divided into $32 = 8\times4$ patches, 
as illustrated in Figure~\ref{fig:QA}(a).
The three dimensional domain is the twisted quarter annulus, decomposed into $128 = 4\times4\times8$ patches as presented in Figure~\ref{fig:QA}(b).

As strongly enforced continuity conditions, we have chosen 
the continuity of the vertex values and the edge averages
for the two dimensional example, and
the continuity of the edge averages
for the three dimensional example.

For the examples with polynomial degree $p=2$, we use a standard MG method based on a hierarchy of nested grids keeping $p$ fixed and use a standard Gauss Seidel (GS) smoother.
For the examples with higher polynomial degree ($p=4$ or $7$), 
we have used 
$p=1$
on all grid levels but the finest grid.
This does not yield nested spaces. 
Thus, we cannot use the canonical embedding and restriction. Instead, we use $L^2$-projections to realize them. On the finest grid, we use a MG smoother suitable
for high-order IgA, namely a variant of the subspace-corrected mass smoother proposed and analyzed in~\cite{HL:HofreitherTakacs:2016a}. For this smoother, it was
shown that a resulting MG method is robust 
with respect to
both the grid size and the polynomial degree. 
However, for $p=1$ or $2$, standard approaches are more efficient.
Thus, we again use this smoother only for the finest level, 
while for all other grid levels we use standard GS smoothers. 
To archive better results,
we have modified the subspace-corrected mass smoother by incorporating a rank-one approximation of the geometry transformation.

For the outer CG or BPCG iteration, we use a zero initial guess, 
and 
the
reduction of the initial residual by the factor $10^{-6}$ 
 as stopping criterion.
The local problems related to the calculation of the $S$-orthogonal basis are solved up to a tolerance of $10^{-12}$. In case of the (MG-MG) version, the local Neumann problems \eqref{equ:localNeumann} in $W_\Delta$ are solved up to a relative error of $10^{-10}$. The number of MG cycles in the preconditioner is fixed. 
For the local Dirichlet problems in the scaled Dirichlet preconditioner, we use $2$ V-cycles. 
The local Neumann problems, 
which appear in the preconditioner of the (MG-MG-S) version, 
are approximately solved by $3$ V-cycles. In the following, we report on the number of CG iterations to solve~\eqref{equ:SchurFinal} and BP-CG iterations for~\eqref{equ:saddlePoint} and the total time, which includes the assembling, the IETI-DP setup and solving phase. 
 
	\begin{table}
	\begin{tabular}{|r||c|r||c|r||c|r||c|r||r||c|r||c|r||c|r||c|r|}\hline
 		  $p=2$   &\multicolumn{2}{c||}{\textbf{D-D} }&\multicolumn{2}{c||}{\textbf{MG-D} }&\multicolumn{2}{|c||}{\textbf{MG-MG} }&\multicolumn{2}{|c||}{\textbf{MG-MG-S}} & $p=7$   &\multicolumn{2}{c||}{\textbf{D-D} }&\multicolumn{2}{c||}{\textbf{MG-D} }&\multicolumn{2}{|c||}{\textbf{MG-MG} }&\multicolumn{2}{|c|}{\textbf{MG-MG-S}}\\ \hline
 		Dofs  & It. & Time & It. & Time & It. & Time & It. & Time & Dofs  & It. & Time & It. & Time & It. & Time & It. & Time \\\hline         \hline
           134421    &  9    &   9.5  &    9    &  7.8    &   9   &   12.5  &  14   &  14.4  &   45753    &   10   &  25.7 &  10    &   26.7   &  10     &   56.7  &   14  &      53.5  \\\hline                                                                                 
           530965    &  10   &   45.4  &   10    &  37.0  &  10   &   54.4  &  15   &  90.1  &   155961   &    11  &   108 &   11   &   110  &    11    &   225   &   15  &      211 \\\hline
           2110485   &   11  &   224  &     11  &    172 &    11 &     272  &   16  &  568    &  572985    &   12   &  498  &  12    &  495  &   12     &   1048  &   17  &     1013  \\\hline
           8415253   &   11  &   1005  &    11   &   762  &   11  &    1181  &  15   &  3394  &  2193465   &    13  &  2384  &   13  &  2265  &    14    & 4427    &   18  &     4344   \\\hline 
          33607701    &   \multicolumn{2 }{c||}{OoM}  &   \multicolumn{2 }{c||}{OoM}    &   13   &      5070 & \multicolumn{2 }{c||}{OoM}  &  8580153   &   \multicolumn{2 }{c||}{OoM}  &  \multicolumn{2 }{c||}{OoM}    &       15  &  18484   &   20  &     19958    \\\hline 
        \end{tabular}
        \caption{Number of outer iterations and timings for the four different formulations using the quarter annulus, see Figure~\ref{fig:QA}(a). GS smoother is used for $p=2$ and $p$-robust subspace corrected mass smoother for $p=7$.
        }
        \label{tab:overview2D}
	\end{table}
	
        \begin{table}
        \begin{tabular}{|r||c|r||c|r||c|r||c|r||r||c|r||c|r||c|r||c|r|}\hline
                  $p=2$   &\multicolumn{2}{c||}{\textbf{D-D} }&\multicolumn{2}{c||}{\textbf{MG-D} }&\multicolumn{2}{|c||}{\textbf{MG-MG} }&\multicolumn{2}{|c||}{\textbf{MG-MG-S}} & $p=4$   &\multicolumn{2}{c||}{\textbf{D-D} }&\multicolumn{2}{c||}{\textbf{MG-D} }&\multicolumn{2}{|c||}{\textbf{MG-MG} }&\multicolumn{2}{|c|}{\textbf{MG-MG-S}}\\ \hline
                Dofs  & It. & Time & It. & Time & It. & Time & It. & Time & Dofs  & It. & Time & It. & Time & It. & Time & It. & Time \\\hline         \hline
            14079     &   11  &   2.6&   11  &      2.5  &   11    &    7.6  & 25    & 7.2  & 40095    &   13 &   29.5 &   13  &    32.8  &  13    &  112 &   23 &    104\\\hline                                                                                 
            86975     &   12  &  19.3&   12  &     19.1  &   12    &   59.1  & 26    & 59.1 & 160863   &   15 &    234 &   15  &    254  &   15    &  659   &  28  &  633 \\\hline
            606015    &   14  &  213 &   14  &    197  &   14     &   484 &    30 &    616  & 849375   &   16 &   2237 &   17  &   2356  &   17    &  5403  &  32 &  5298 \\\hline
            4513343   &   \multicolumn{2 }{c||}{OoM}  &   16  &   2764  &   16    &   5244 &    35 &  11657 & 5390559&   \multicolumn{2 }{c||}{OoM}  &      \multicolumn{2 }{c||}{OoM}      &    19     &        45243 &  37 &    52831 \\\hline 
        \end{tabular}
        \caption{Number of outer iterations and timings for the four different formulations using the twisted quarter annulus, see Figure~\ref{fig:QA}(a). GS smoother is used for $p=2$ and $p$-robust subspace corrected mass smoother for $p=4$. }
        \label{tab:overview3D}
        \end{table}
	
	The algorithm is realized with the open source C++ library G+Smo \cite{gismoweb}, which uses the linear algebra facilities of the Eigen library \cite{HL:eigenweb}. We utilize the PARDISO 5.0.0 Solver \cite{HL:PARDISO500} for performing the LU factorizations.\footnote{Our code is compiled with the \texttt{gcc 4.8.3} compiler with optimization flag \texttt{-O3}. The results are obtain on the RADON1 cluster at Linz. We use a single core of a node, equipped with 2x Xeon E5-2630v3 ``Haswell'' CPU (8 Cores, 2.4Ghz, 20MB Cache) and 128 GB RAM.}

  In Table~\ref{tab:overview2D}, we summarize the results for the two dimensional domain for $p=2$ and $7$.
   We observe that replacing the direct solver in the preconditioner with two MG V-cycles does not change the number of outer iterations. Moreover, going from the Schur complement  to the saddle point formulation and using BPCG there, leads only to a minor increase in the number of outer iterations. In all cases, the logarithmic dependence of the condition number on $h$ is preserved. The advantage of the formulation using only MG, especially (MG-MG), is its smaller memory footprint, therefore, the possibility of solving larger systems. However, the setting with the best performance is (MG-D). Concluding, for small polynomial degrees and using the GS smoother, (MG-MG) gives reasonable trade off between performance and memory usage and for larger polynomial degrees, this setting can be still recommended if memory consumption is an issue.

   In the case $p=2$, for the inner iterations, we have observed
   that the CG needed on average 8 iterations to compute $\widetilde{g}$, the calculation of the $S$-orthogonal basis needed on average 14 iterations and the solution of~\eqref{equ:localNeumann} 
  required on average 10 iterations. For the second case, $p=7$, we needed 9 iterations to compute $\widetilde{g}$, 13 iterations for the calculation of the $S$-orthogonal basis and 10 iterations for the solutions of~\eqref{equ:localNeumann}.
   Here and in what follows, we have taken the average over the patches, the individual levels and the individual steps of the outer iteration. 
    We mention
   that the number of inner iterations was only varying slightly.

   In Table~\ref{tab:overview3D}, we summarize the results for the three dimensional domain and for $p=2$ and $4$.
   We observe that replacing the direct solver in the preconditioner with two MG V-cycles does not change the number of outer iterations. 
   We further observe that the results behave similar to the one of the two dimensional case. However, the number of iterations almost doubled when using BPCG for (MG-MG-S). 
   In all cases, the logarithmic dependence of the condition number on $h$ is preserved.
   The advantage of the formulation using only MG, especially (MG-MG), is its smaller memory footprint, therefore the possibility of solving larger systems.
   The best performance is obtained sometimes by (D-D) and sometimes by (MG-D), where both approaches are comparable.

   Concerning the inner iterations, for $p=2$, we need on average 15 CG iterations to compute~$\widetilde{g}$, 22 CG iterations to build up each $S$-orthogonal basis function, and 18 CG iterations to solve \eqref{equ:localNeumann}. In the case of $p=4$, we needed on average only 10 iterations to compute~$\widetilde{g}$, 14 iterations for the construction of the $S$-orthogonal basis functions, and 11 iterations for solving \eqref{equ:localNeumann}.

\section*{Acknowledgments}
This work was supported by the Austrian Science Fund (FWF) under the grant W1214, project DK4. This support is gratefully acknowledged. 

\bibliographystyle{abbrv}

\bibliography{DD24_proceedings.bib}

\end{document}